\newcommand{\la }{\langle}
\newcommand{\ra }{\rangle}

\newcommand{\Exp}{{\rm Exp}}

\newcommand{\mcc}{{\mathcal C}}

\newcommand{\mcf}{{\mathcal F}}

\newcommand{\mcfs}{{\mathcal F}_s}

\newcommand{\mch}{{\mathcal H}}

\newcommand{\mbc}{{\mathbb C}}

\newcommand{\mbe}{{\mathbb E}}

\newcommand{\mbr}{{\mathbb R}}

\documentclass[10pt,draft,reqno]{amsart}
     \makeatletter
     \def\section{\@startsection{section}{1}%
     \z@{.7\linespacing\@plus\linespacing}{.5\linespacing}%
     {\bfseries
     \centering
     }}
     \def\@secnumfont{\bfseries}
     \makeatother
\newtheorem{theorem}{Theorem}[section]

\newtheorem{proposition}[theorem]{Proposition}

\theoremstyle{definition}

\theoremstyle{remark}

\numberwithin{equation}{section}
\begin{document}

\title[Gaussian Radon for Wiener Space]{The Gaussian Radon Transform in Classical Wiener Space}

\author{Irina Holmes}
\address{School of Mathematics, Georgia Institute of Technology, 686 Cherry Street
Atlanta, GA 30332-0160, USA}
\email{irina.c.holmes@gmail.com }

\author{Ambar N.~Sengupta}
 \thanks{* This research is supported in part by NSA grant
  H98230-13-1-0210.}
\address{A. N. Sengupta: Department of Mathematics, Louisiana State University, Baton Rouge, LA 70810, USA}
\email{ambarnsg@gmail.com}



\begin{abstract}
We  study the Gaussian Radon transform in the classical Wiener space of Brownian motion. We determine explicit formulas for transforms of Brownian functionals specified by stochastic integrals. A Fock space decomposition is also established for  Gaussian measure conditioned to  closed affine subspaces in Hilbert spaces.
\end{abstract}

\maketitle


 \section{Introduction}\label{s:intro}

 The Gaussian Radon transform $G\phi$ of a function $\phi$ on a linear space associates to each closed affine  subspace $L$ the integral of $\phi$ with respect to a Gaussian measure on $L$. In this paper we study this transform for the classical Wiener space, which is the Banach space of continuous paths starting at a given basepoint. We  also establish a Fock space decomposition in this context, expanding suitable functions $\phi$ in terms of symmetric  tensors.  
 
 The classical Radon transform \cite{Radon} associates to a function $\phi$ on Euclidean space $\mbr^n$ its integrals over hyperplanes. An infinite dimensional analog of this transform  was introduced in  Mihai and Sengupta\cite{MS} (an early study of this type was carried out by Hertle \cite{Hert}). Further development of this theory, including the establishment of counterparts of the classical Helgason support theorem \cite{H}, has been carried out in the works of Becnel \cite{BecGR2010, BecSen12}, Bogachev and Lukintsova \cite{BogLukin,   Lukin}, and by us  \cite{HolSen2012}.  A  comprehensive account of Gaussian measures on infinite dimensional spaces is presented in Bogachev  \cite{Bog}. The  theory of abstract Wiener spaces, originating in the work of Gross \cite{Gr}, is presented in the monograph of Kuo \cite{Ku1} as well as in many other works.
 
 The present paper is organized as follows. We begin in section \ref{s:grt} with a self-contained description of the Gaussian Radon transform in Banach spaces and introduce necessary notation and concepts.  In infinite dimensions it is Gaussian measure $\mu$ that needs to be used as there is no useful notion of Lebesgue measure. An important  issue  is that $G\phi$ depends on the pointwise defined function $\phi$ and not simply its $\mu$-almost-everywhere specified representative. For any  bounded continuous function $\phi$ the transform $G\phi$ is meaningful without further specifications. In section \ref{s:grt}  we also present a more detailed discussion of the case of affine subspaces of finite codimension. In section \ref{s:fock} we present a Fock space decomposition for functions with respect to Gaussian measure on closed affine subspaces.  We turn next in section \ref{s:cw} to a description of the classical Wiener space, formulated in a manner useful for our objectives, along with the Gaussian Radon transform in this setting. In section \ref{s:cw} we apply   the machinery developed to work out Gaussian Radon transforms of some specific functions on the classical Wiener space, such as those given by iterated It\^o integrals.

 \section{The Gaussian Radon Transform}\label{s:grt}
 
 In this section we present a description of the Gaussian Radon transform, referring to \cite{HolSen2012}   for further details and proofs. We begin with a   review of  the framework of an abstract Wiener space.
 
\subsection{The Gaussian measure $\mu$}   We work with a real separable Banach space $B$, with norm denoted $|\cdot|$, and a Borel probability measure $\mu$ on $B$ for which the distribution of  every continuous linear functional $x^*:B\to\mbr$ is Gaussian with mean $0$, the variance being positive if $x^*$ is not zero.  The characteristic function of $x^*$ is therefore given by
\begin{equation}\label{E:chxstar}
\int_Be^{itx^*}\,d\mu=e^{-\frac{t^2}{2}|\!|x^*|\!|_{L^2(\mu)}^2},
\end{equation} 
for all $t\in\mbr$, 
because the variance of $x^*$ is $\mbe_{\mu}[{x^*}^2]$, which is the $L^2$-norm squared. By a  Cameron-Martin space for $(B, \mu)$  we shall mean  a Hilbert space $H$ along with a linear isometry
$$I:H\to L^2(B,\mu)$$
whose image contains $B^*$, the set of all continuous linear functionals on $B$, as a dense subspace. For example, $H$ could be   the closure of $B^*$ in $L^2(B, \mu)$, with $I$ being the inclusion map.  We   use the notation
\begin{equation}\label{E:defHBstar}
H_{B^*}=\{h\in H\,:\, I(h)\in B^*\}.
\end{equation}
By assumption, the image $I(H_{B^*})=B^*$ is dense in the range of the linear isometry $I$, and so $H_{B^*}$ is dense in $H$.
 For $h\in H$ we can choose a sequence $x_n^*\in B^*$ converging to $I(h)$ in $L^2(\mu)$-norm; the charactertistic function of $x_n^*$ is given by (\ref{E:chxstar}) with $x^*=x_n^*$, and then by using dominated convergence we have
\begin{equation}\label{E:chIhstar}
\int_Be^{itI(h)}\,d\mu=e^{-\frac{t^2}{2}|\!|I(h)|\!|_{L^2(\mu)}^2},
\end{equation} 
for all $t\in\mbr$. Thus  $I(h)$ is Gaussian with mean $0$ and variance $|\!|h|\!|^2$, for every $h\in H$.   
There is also a continuous linear injective map
\begin{equation}\label{E:defiHB}
i:H\to B: h\mapsto \int_B (Ih)(\omega) \omega\,d\mu(\omega),
\end{equation}
where the integrand is a $B$-valued function  and the  integral is a  Bochner integral.   In particular
\begin{equation}\label{E:iHBhxstar}
\la x ^*, i(h)\ra= \int_B Ih(\omega) x^*(\omega) \,d\mu(\omega),
\end{equation}
for all $x^*\in B^*$ and $h\in H$. Taking $x^*=I(h_1)$, where $h_1\in H_{B^*}$, we obtain
\begin{equation}\label{E:iHBhh_1}
\la I(h_1), i(h)\ra= \la I(h_1), I(h)\ra_{L^2(\mu)}=\la h_1, h\ra_{H}
\end{equation}
for all $h\in H$, where on the left we have the evaluation of $I(h_1)\in B^*$ on $i(h)$ and in the second expression on the right we have stressed through the subscript that the inner product is in $H$.  In particular, if $i(h)=0$ then $h$ is orthogonal to  every $h_1\in H_{B^*}$, and so, since $H_{B^*}$ is dense in $H$,  $h=0$. Thus  the mapping {\em   $i$ is injective.}

  Another perspective on $H_{B^*}$ is also useful: an element $h_1\in H$ lies in the subspace $H_{B^*}$ if and only if  the mapping
  $$H\to\mbr: h\mapsto \la h,h_1\ra$$
  is continuous with respect to the norm $|\cdot|$ on $H$ obtained by the identification of $H$ with the subspace $i(H)$ of $B$.
  
  The system $i:H\to B$, along with the Gaussian measure $\mu$ on $B$,  is called an {\em abstract Wiener space}. This notion was introduced by Gross \cite{Gr}.

\subsection{The measure on affine subspaces} A closed affine subspace of $H$ is of the form
$$L=p+L_0,$$
where $L_0$ is a closed subspace of $H$ and $p$ a point $H$.  In  \cite[Theorem 2.1]{HolSen2012}  we have constructed, for every closed affine subspace $L$ of $H$,  a Borel probability measure $\mu_L$ on $B$  such that for any $h\in H_{B^*}$
\begin{equation}\label{E:muLdef1}
\int_Be^{iI(h)}\,d\mu_L=e^{i\la h_L,h\ra -\frac{1}{2}|\!|P_{L_0}h|\!|^2}
\end{equation}
where   $h_L$ is the point on $L$ closest to $0\in H$ and
\begin{equation}\label{E:PL0}
P_{L_0}: H\to H
\end{equation}
is the orthogonal projection onto the closed subspace 
\begin{equation}\label{E:L0LhL}
L_0=-h_L+L.
\end{equation}
 The measure $\mu_L$ is concentrated on the closure $\overline{L}$ of $L$ in the Banach space $B$:
\begin{equation}\label{E:muLL1}
\mu_L({\overline L})=1.
\end{equation}
(This is also proved in \cite[Theorem 2.1]{HolSen2012}.) One other fact to note here is the nature of the closure $\overline{L}$.  Suppose $L_0=\ker f$ where $f\in H^*$ is the restriction to $H$ of an element $f_B\in B^*$;  then 
\begin{equation}\label{E:kerfB}
\overline{L}_0=\ker f_B.
\end{equation} 
If $f$ is not continuous with respect to the norm on $B$ then  $\overline{L}=B$. This is \cite[Proposition 5.2]{HolSen2012}.

\subsection{The mapping $I_L$} 

Let $L$ be a closed affine subspace of $H$. The characteristic function of a Gaussian variable $X$ is
  $$\mbe[e^{itX}] =e^{it\mbe[X]-\frac{t^2}{2}{\rm var}(X)},
  $$
  where ${\rm var}(X)$ is its variance.  Then from (\ref{E:muLdef1}), with $th$ in place of $h$, we see that when $h\in H_{B^*}$ the random variable $I(h)$ has Gaussian distribution with respect to $\mu_L$, having
 mean and variance
\begin{equation}\label{E:meanvarIh}
\mbe_{\mu_L}\Bigl(I(h)\Bigr)=\la h_L, h\ra\qquad\hbox{and}\qquad {\rm Var}_{\mu_L}\Bigl(I(h)\Bigr)=|\!|P_{L_0}h|\!|^2.
\end{equation}
 Consider now the mapping
 \begin{equation}\label{E:ImuL}
 I^0_L:H_{B^*}\to L^2(B,\mu_L): h\mapsto I(h).
 \end{equation}
 Although it appears to be just the mapping $I$, there is in fact  a significant difference because here $I(h)$ is taken as an element of $L^2(\mu_L)$ (thus, for example, $I_L(h)$ might be $0$ even if $I(h)$ is not $0$). Clearly, $I^0_L$ 
 is linear, and, moreover, 
 \begin{equation}
 |\!| I^0_L(h)|\!|^2=\mbe_{\mu_L}\bigl(I^0_L(h)^2\bigr)=  {\rm Var}_{\mu_L}\bigl(I(h)\bigr)+ \Bigl(\mbe_{\mu_L}\bigl(I(h)\bigr)\Bigr)^2\leq (1+|\!|h_L|\!|^2)|\!|h|\!|^2 
 \end{equation}
 for all $h\in H_{B^*}$.
 Thus   $I^0_L$, being    continuous  with respect to the norm  of $H$,   extends to a continuous linear mapping
 \begin{equation}\label{E:defILh}
 I_L:H \to L^2(B,\mu_L).
 \end{equation}
 In an intuitive sense, $I_L(h)$ is the ``restriction'' of $I(h)$ to the closed affine subspace $\overline{L}$, the closure of $L$ in $B$.    Let us note again, however,  that for a general element $h\in H$ the random variable
 $$I(h)$$
 is an element if $L^2(\mu)$, thus defined almost everywhere with respect to $\mu$, and cannot generally be ``restricted'' to $\overline{L}$.  For $h\in H$, the random variable $I_L(h)$ has  characteristic function given by
\begin{equation}\label{E:muLdef}
\int_Be^{iI_L(h)}\,d\mu_L=e^{i\la h_L,h\ra -\frac{1}{2}|\!|P_{L_0}h|\!|^2}, 
\end{equation}
with $h_L$ being as before, 
as can be seen by applying (\ref{E:muLdef1}) to a sequence of elements $h_n\in H_{B^*}$ converging in $H$ to $h$ and evaluating the limit using the dominated convergence theorem. 

As with (\ref{E:meanvarIh}), we see that  
 $I_L(h)$ is Gaussian and  
\begin{equation}\label{E:meanvarILh2}
\mbe_{\mu_L}\bigl(I_L(h)\bigr)=\la h_L, h\ra\qquad\hbox{and}\qquad {\rm Var}_{\mu_L}\bigl(I_L(h)\bigr)=|\!|P_{L_0}h|\!|^2.
\end{equation}
(The difference between this and  (\ref{E:meanvarIh}) is that now $h$ is any element of $H$, not restricted to being in $H_{B^*}$.) 
In particular, {\em if $h$ is orthgonal to $L$ (more precisely, to the subspace $L-h_L$) then $I_L(h)$ has zero variance and is thus just the constant $\la h,h_L\ra$. }

From the variance formula in (\ref{E:meanvarILh2}) with $h=th_1+sh_2$, where $h_1, h_2\in H$ and $t,s\in\mbr$, it follows that the covariance of $I_L(h_1)$ and $I_L(h_2)$ is given by
\begin{equation}\label{E:covILH12}
{\rm Cov}\bigl(I_L(h_1), I_L(h_2)\bigr)= \la P_{L_0}h_1, P_{L_0} h_2\ra.
\end{equation}

As with any Gaussian we have
\begin{equation}
\label{E:Gaussfundaform}
\int_B e^{z I_L(h)}\,d\mu =e^{z\la h_L,h\ra +\frac{z^2}{2}|\!|P_{L_0}h|\!|^2}, 
\end{equation}
for all $h\in H$ and all $z\in\mbc$.

\subsection{The Gaussian Radon transform $G\phi$}
Consider a measurable function $\phi$ on $B$. The Gaussian Radon transform $G\phi$ is a function  that associates to a closed affine subspace $L\subset H$,  the expectation of $\phi$ with respect to the measure $\mu_L$:
\begin{equation}\label{E:defGFL}
G\phi|_L\stackrel{\rm def}{=}G\phi (L)\stackrel{\rm def}{=}\int_B \phi\,d\mu_L.
\end{equation}
 Thus the domain of $G\phi$ is the set of all $L$ for which this integral is finite.  If $\phi$ is a bounded Borel function then $G\phi$ is defined on  the set of all closed affine subspaces.  
 
 We observe that  $G$ is defined for a  pointwise specified function, and not on $L^2(\mu)$. For example, a function $\psi$ might be $0$ almost everywhere with respect to $\mu$ but $G\psi$ need not be $0$.
 
 As a special case consider $\phi=I(h)$ where $h\in H_{B^*}$; let us note that this means that $I(h)$ is a continuous function on $B$.  From the formula (\ref{E:meanvarILh2}) we conclude that the Gaussian Radon transform of $I(h)$ is
 \begin{equation}\label{E:GRTIh}
 GI(h)|_L=\la h_L,h\ra
 \end{equation}
 where $h_L$ is the point on the closed affine subspace $L$ that is of smallest norm.
 
 \subsection{Conditioning a function to an affine subspace}   As we have seen above, for each vector $h\in H$ there is a function $I(h)\in L^2(\mu)$ and also a function $I_L(h)\in L^2(\mu_L)$. If $I_L(h)=0$ then $I_L(th)=0$ for all $t\in\mbr$ and so by (\ref{E:muLdef}),   $P_{L_0}h=0$ and $\la h_L,h\ra=0$, and so $h$ is a vector in $L_0^\perp$ that is orthogonal to $h_L$.    On the other hand, if $I(h)=0$ then $h=0$ (because $I:H\to L^2(\mu)$  is an isometry) and so $I_L(h)=0$. Thus there is a continuous linear map
 \begin{equation}
I_L\circ I^{-1}:L^2(\mu)_0\to L^2(\mu_L): I(h)\mapsto I_L(h)
\end{equation}
where $L^2(\mu)_0$ is the image of $I$ inside $L^2(\mu)$. This provides a way of obtaining functions in $L^2(\mu_L)$ from functions in $L^2(\mu)_0$. 

\subsection{The point $h_L$ for finite codimension subspaces}  Consider a closed affine subspace of $H$ that is specified as a level surface of a finite number of linear functionals; specifically, consider the affine subspace $L_c\subset H$ given by
\begin{equation}\label{E:Lc}
L_c=\{h\in H\,:\, \la v_1,h\ra=c_1,\ldots,  \la v_m, h\ra=c_m\},
\end{equation}
where $v_1,\ldots, v_N\in H$ and $c=(c_1,\ldots, c_m)\in \mbr^m$.  
(Such subspaces were used in the 
context of machine learning in our work  \cite{HolSen2014}, where formulas of the type discussed below in Proposition \ref{P:hLforcodimfin} were obtained.)  Then we have the following result in linear algebra.

\begin{proposition}\label{P:hLforcodimfin} Let $H$ be a real Hilbert space
and $F:H\to\mbr^m$ a linear map given by
\begin{equation}
F(h)=\Bigl(\la v_1,h\ra,\ldots, \la v_m, h\ra\Bigr)
\end{equation}
where $v_1,\ldots v_n$ are linearly independent vectors in $ H$  and $c=(c_1,\ldots, c_m)\in \mbr^m$. Then the mapping $FF^*:\mbr^m\to\mbr^m$ is invertible and  the point on the closed affine subspace $L_c=F^{-1}(c)$ closest to the origin $0$ is
\begin{equation}\label{E:hLF}
h_{L_c}=F^*(FF^*)^{-1}(c).
\end{equation}
The matrix for $FF^*$ is
\begin{equation}\label{E:matrFFstar}
FF^*=\left[\begin{matrix} \la v_1, v_1\ra & \ldots & \la v_m, v_1\ra \\
\vdots &\vdots &\vdots  \\
\la v_1, v_m\ra &\ldots & \la v_m, v_m\ra \end{matrix}\right],\end{equation}
and  the point $h_{L_c}$ can be expressed as
\begin{equation}\label{E:hLF2}
h_{L_c} =\sum_{k=1}^m a_kv_k
\end{equation}
 where 
 \begin{equation}\label{E:ak}
a_k  = \sum_{j=1}^m  [(FF^*)^{-1}]_{kj}c_j.
\end{equation}
If the vectors $v_j$ are orthonormal then 
\begin{equation}\label{E:matrFFstar2}
h_{L_c}=F^*(c)=\sum_{k=1}^mc_kv_k.
\end{equation}
\end{proposition}
\noindent\underline{Proof}.    For any $b\in\mbr^m$ and $w\in H$ we have:
\begin{equation}
 \la b, F(w)\ra_{\mbr^m}=\sum_{k=1}^m b_k\la v_k,w\ra =\left\la\sum_{k=1}^mb_kv_k, w\right\ra
\end{equation}
Hence the adjoint $F^*:\mbr^m\to H$  is given by
\begin{equation}\label{E:Fstarb}
F^*(b)=\sum_{k=1}^mb_kv_k,
\end{equation}
for all $b\in\mbr^m$.  In particular, taking $b=e_j$, the $j$-th standard basis vector in $\mbr^m$ we have
\begin{equation}\label{E:Fstarb2}
F^*(e_j)=v_j.
\end{equation}
Since $v_1,\ldots, v_m$ are linearly independent we see that $\ker F^*=\{0\}$ which implies also that $FF^*$ is invertible because any vector $b$ that lies in $\ker (FF^*)$ satisfies 
$$|\!|F^*b|\!|^2=\la  b, FF^*b\ra =0.$$
Thus the   definition of the vector $h_{L_c}$ given in (\ref{E:hLF}) is meaningful. The definition shows directly that $h_{L_c}$ lies in $F^{-1}(c)$. Furthermore,  if $w$ is any point on $F^{-1}(c)$  then 
$$\la w,h_{L_c}\ra =\la F(w), (FF^*)^{-1}(c)\ra_{\mbr^m}=\la c, (FF^*)^{-1}(c)\ra,$$
which, being independent of the choice of $w$ is equal to the value obtained on taking $w$ to be the point $h_{L_c}$; consequently
$$\la w-h_{L_c}, h_{L_c}\ra=0.$$
Thus $h_{L_c}$ is a point on $L_c$ that is orthogonal to all vectors in the subspace $L_c-h_{L_c}$. This implies that the vector $h_{L_c}$ forms the perpendicular from the origin onto the affine subspace $L_c$. Hence $h_{L_c}$ is the point on $F^{-1}(c)$ of smallest norm (which may also be seen by noting that $|\!|h_{L_c}|\!|^2=\la w, h_{L_c}\ra\leq    |\!|w|\!||\!|h_{L_c}|\!|$ for all $w\in L_c$, from which it follows that $|\!|h_{L_c}|\!|\leq |\!|w|\!|$ for such $w$.)

To express things in coordinates we apply $F$ to $F^*(b)$ as given in (\ref{E:Fstarb})  to obtain
\begin{equation}\label{E:FFstb}
FF^*(b)=\sum_{k=1}^mb_k\left(\sum_{j=1}^m\la v_k,v_j\ra e_j\right)
=\sum_{j=1}^m\left(\sum_{k=1}^m \la v_k,v_j\ra b_k\right)e_j.
\end{equation}
Thus the $(j,k)$-th entry of the matrix for $FF^*$ is $\la v_k,v_j\ra$.
Using (\ref{E:hLF}) we have
\begin{equation}\label{E:hLF3}
h_{L_c}=F^*(FF^*)^{-1}(c)=\sum_{j,k=1}^m[(FF^*)^{-1}]_{kj}c_jF^*(e_k)
\end{equation}
and this establishes the expression (\ref{E:hLF2}) for $h_{L_c}$ upon using $F^*(e_k)=v_k$ as observed in (\ref{E:Fstarb2}).

Finally, if the vectors $v_j$ are orthonormal then $FF^*=I$ and then the expression for $h_{L_c}$ reduces to $F^*(c)$. \fbox{QED}

The explicit formula for $h_{L_c}$ given above makes it possible for us to write down a formula for the Gaussian Radon transform of $I(h)$ for $h\in H_{B^*}$:
\begin{equation}\label{E:GIhhLc}
\begin{split}
GI(h)|_{L(c)} &=\la h_{L_c}, h\ra \\
&=\la (FF^*)^{-1}(c), F(h)\ra_{\mbr^m}\\
& = \sum_{j,k=1}^m[(FF^*)^{-1}]_{jk}c_k\la  v_j, h\ra.  
\end{split}
\end{equation}
Here we have used the formula (\ref{E:GRTIh}) for $GI(h)|_L$ as well as the expression for $h_{L_c}$ from (\ref{E:matrFFstar2}). In the special case where the vectors $v_j$ are orthonormal this simplifies further to
\begin{equation}\label{E:GIhhLcorth}
GI(h)|_{L(c)}=\sum_{k=1}^mc_k\la  v_k, h\ra
\end{equation}

\section{Fock space decomposition}\label{s:fock}

There is a well known unitary isomorphism between $L^2$ of a Gaussian space and a completed tensor algebra. In this section we construct such an isomorphism for the conditioned measure $\mu_L$. We work with a centered Gaussian measure $\mu$ on a separable real Banach space $B$, as  discussed before,  and $H$ is a Cameron-Martin Hilbert space for $(B, \mu)$.

\subsection{Tensor powers} The algebraic tensor power
$$\otimes_{\rm alg}^nH$$
which is a vector space spanned by vectors of the form $v_1\otimes\ldots\otimes v_n$ with $v_1,\ldots, v_n\in H$; in slightly greater details, if $\{b_\alpha\}_{\alpha\in I}$ is a vector space basis of $H$ then the abstract expressions $b_{\alpha_1}\otimes\ldots\otimes b_{\alpha_n}$, with $\alpha_1,\ldots,\alpha_n$ running over $I$, form a basis of the vector space $\otimes_{\rm alg}^nH$. On this vector space there is a unique inner product specified by the requirement that
\begin{equation}
\la v_1\otimes\ldots\otimes v_n, w_1\otimes\ldots\otimes w_n\ra_n=\prod_{j=1}^n\la v_j, w_j\ra.
\end{equation}
The completion of $\otimes_{\rm alg}^nH$ with respect to this inner product is the Hilbert space
$$\otimes^nH,$$
which is the n-th tensor power of $H$ in the Hilbert-space sense.  Let $S_n$ denote the symmetric group $S_n$ of all permutations on $\{1,\ldots, n\}$. Then for each $\sigma\in S_n$ there is a  unitary isomorphism $U_{\sigma}$ on $H^{\otimes n}$ specified by:
$$U_{\sigma}:H^{\otimes n}\to H^{\otimes n}: v_1\otimes\ldots \otimes v_n\mapsto v_{\sigma^{-1}(1)}\otimes\ldots \otimes v_{\sigma^{-1}(n)}.$$
The set of all $x\in H^{\otimes n}$ that satisfy
$$U_{\sigma}x=x\qquad\hbox{for all $\sigma\in S_n$}$$
is a closed subspace of $H^{\otimes n}$ called the $n$-th {\em symmetric tensor power} of $H$; it is denoted
$$H^{\hat\otimes n}.$$
In the case $n=0$ we take this space to be just $\mbr$ with its usual inner-product.  

\subsection{Symmetric Fock space}
Next we consider the {\em symmetric tensor algebra}, which is the vector space direct  sum
\begin{equation}
S(H)=\Sigma_{n\geq 0} H^{\hat\otimes n}.
\end{equation}
Thus each element of $S(H)$ is of the form
$$\sum_{n\geq 0}x_n,$$
where $x_n\in H^{\hat\otimes n}$ and all except finitely many $x_n$ are $0$.

On $S(H)$ there are different natural choices for an inner product; we work with the inner product specified by
\begin{equation}
\left\la\sum_{n\geq 0}x_n,\sum_{m\geq 0}y_m\right\ra=\sum_{n\geq 0}n!\la x_n, y_n\ra_n,
\end{equation} 
wherein $x_n, y_n\in H^{\hat\otimes n}$. The completion of $S(H)$ with respect to this inner product is a Hilbert space called the   {\em symmetric Fock space} over $H$ and we will denote it by
$${\mathcal F}_s(H).$$
Finally, we can consider the complexification of this, which is a complex Hilbert space that we denote by
$${\mathcal F}_s(H)_c.$$
\subsection{Exponential vectors} Inside this Hilbert space are certain useful vectors called {\em exponential vectors} or {\em coherent state} vectors; these are of the form
\begin{equation}\label{E:Exp}
\Exp(v)=\sum_{n=0}^\infty\frac{1}{n!}v^{\otimes n},
\end{equation}
for any $v\in H$. We observe that
\begin{equation}\label{E:Expvwip}
\left\la\Exp(v),\Exp(w)\right\ra=\sum_{n=0}^\infty n!\frac{1}{n!^2}\la  v^{\otimes n}, w^{\otimes n}\ra_n=e^{\la v,w\ra}
\end{equation}
for all $v, w\in H$. The other essential fact is that the linear span of the vectors $\Exp(v)$ is dense in $\mcfs(H)$.

\subsection{Fock space decomposition}
Now let us bring in the Banach space $B$ with  centered Gaussian measure $\mu$, for which $H$ is a Cameron-Martin space.
There is a fundamental unitary isomorphism
\begin{equation}\label{E:fockspaceu}
U:\mcfs(H)_c\to L^2(B,\mu)
\end{equation}
which is uniquely specified by the requirement that it map $\Exp(v)$ to a multiple of the exponential function $e^{I(v)}$; more specifically
\begin{equation}\label{E:IvExp}
U\Bigl(\Exp(h)\Bigr)= e^{I(h)-\frac{1}{2}|\!|h|\!|^2}\qquad\hbox{for all $h\in H$. }\end{equation}
The surjectivity of $U$  follows from the fact that the random variables of the form $e^{I(h)}$ span a dense subspace of $L^2(B, \mu)$. The relation (\ref{E:Expvwip}) implies that $\Exp$ is an isometry.

\subsection{Fock space decomposition for $\mu_L$}

Now let $L$ be a closed affine subspace of $H$.  We denote by $h_L$ the point on $L$ closest to the origin
\begin{equation}
|\!|h_L|\!|=\inf\{|\!|v|\!|\,:v\in L\}.
\end{equation}
The vector $h_L$ is orthogonal to $L$; more precisely,
\begin{equation}
h_L\in L_0^\perp.
\end{equation}

We have then the following result on Fock space decomposition of $L^2(B,\mu_L)$:

\begin{theorem}\label{T:FockBmuL} Let $H$ be a Cameron-Martin Hilbert space associated with the separable real Banach space $B$ equipped with a centered Gaussian measure $\mu$. Let $L$ be a closed affine subspace of $H$,  $h_L$ the point on $L$ closest to the origin, and $L_0$ the closed subspace of $H$ given by
\begin{equation}
L_0=L-h_L.  
\end{equation}
Let $\mu_L$ be the probability measure on $B$ given by  (\ref{E:muLdef1}) and $I_L:H\to L^2(B,\mu_L)$ the continuous linear map given in  (\ref{E:defILh}).
Then there is a unique unitary isomorphism
\begin{equation}
U_L:\mcfs(L_0)_c\to L^2(B,\mu_L)
\end{equation}
for which
\begin{equation}\label{E:ULspecif}
U_L\Bigl(\Exp(v)\Bigr)= e^{I_L(v) -\frac{1}{2}|\!|v|\!|^2}
\qquad\hbox{for all $v\in L_0$.}
\end{equation}
\end{theorem}
Taking $v=0$ we see that $U_L$ maps the `tensor' $1\in \mbc\subset \mcfs(L_0)_c$ to the constant function $1\in L^2(B,\mu_L)$. Moreover, replacing $v$ by $tv$ and `comparing the coefficient' of $t$ on both sides we have
\begin{equation}\label{E:ULIL}
U_L(v)=I_L(v)\qquad\hbox{for all $v\in L_0$,}
\end{equation}
showing that $U_L$ reduces to the mapping $I_L$ on $L_0\subset\mcfs(L_0)_c$.

\noindent\underline{Proof}.  Let us verify that $U_L$ preserves inner-products for the coherent vectors $\Exp(v)$.  For any $v$ and $w$ in the closed subspace $L_0\subset H$ we have
\begin{equation}\label{E:ULip}
\begin{split}
&\int_B e^{I_L(v) -\frac{1}{2}|\!|v  |\!|^2}e^{I_L(w) -\frac{1}{2}|\!|w |\!|^2}\,d\mu_L \\
&=e^{ -\frac{1}{2}\Bigl(|\!|v |\!|^2+|\!|w |\!|^2\Bigr)}\int_B e^{I_L(v+w)}\,d\mu_L\\
&=e^{-\frac{1}{2}\bigl(|\!|v |\!|^2+|\!|w |\!|^2\bigr) +\la v+w,h_L\ra+\frac{1}{2}|\!|P_{L_0}(v+w)|\!|^2}\\
&\qquad\hbox{(using the Gaussian formula (\ref{E:Gaussfundaform}))}\\
&=e^{-\frac{1}{2}\bigl(|\!|v |\!|^2+|\!|w |\!|^2\bigr)  +\frac{1}{2}|\!| v+w|\!|^2}\\
&\qquad\hbox{(since $v, w\in L_0$, and $h_L\in L_0^\perp$; also   $P_{L_0}v=v$ and $P_{L_0}w=w$)}\\
&=e^{\la v,w\ra},
\end{split}
\end{equation}
and this is equal to the inner product
$$\la\Exp(v),\Exp(w)\ra=e^{\la v,w\ra}.$$
This shows that $U_L$, specified on the vectors of the form $\Exp(v)$ by 
(\ref{E:ULspecif}), is well-defined as an isometric linear mapping on the complex linear span of $\{\Exp(v)\,:v\in H\}$ in $\mcfs(H)_c$. Thus, since this subspace is dense in $\mcfs(H)_c$, the mapping $U_L$ extends uniquely to a continuous linear mapping  $\mcfs(H)_c\to L^2(B, \mu_L)$. The functions $e^{I_L(v)}$ span a dense subspace in $L^2(B,\mu_L)$ and so $U_L$ maps $\mcfs(H)_c$ onto $L^2(B, \mu_L)$. \fbox{QED}

 \section{The classical Wiener space and the conditioned measure}\label{s:cw}
 
 In this section we review facts about the classical Wiener space and then extend some of the considerations to the measure $\mu_L$. For more details and further development of the theory we refer to the works of Gross \cite{Gr} and  Kuo \cite{Ku1}.
 
 \subsection{The classical Wiener space} The Banach space here is  
$$\mcc_0=C_0[0,1],$$
 the space of continuous functions $\omega: [0,1]\to\mbr$, with initial value $\omega(0)=0$, with the norm being the supremum norm, and the measure $\mu$ being standard Wiener  measure. In this section we shall review standard facts with a view to use in the contex of the Gaussian Radon transform. For each $t\in [0,1]$ we have  the evaluation map:
\begin{equation}\label{E:defBt}
B_t:\mcc_0\to\mbr:\omega\mapsto B_t(\omega)=B(t;\omega)=\omega(t).
\end{equation}
 The process
  $$t\mapsto B_t$$
   is standard Brownian motion under the measure $\mu$.  This means that each $B_t$ is a Gaussian random variable of mean $0$ and variance $t$, and for any $t_1,\ldots, t_n\in [0,1]$ the variable $(B_{t_1},\ldots, B_{t_n})$ is an $\mbr^n$-valued Gaussian, and
 \begin{equation}\label{E:BsBt}
 \mbe_{\mu}\Bigl(B_sB_t\Bigr)=\min\{s,t\} \qquad\hbox{for all $s,t\in [0,1]$.} \end{equation}

\subsection{Reproducing kernel} A convenient choice for the Cameron-Martin space $\mch_0$  is given as follows. For each $s\in [0,1]$ let  $K_s$ be the function on $[0,1]$ given by
 \begin{equation}\label{E:hsminst}
K_s:[0,1]\to\mbr: t\mapsto  K_s(t)\stackrel{\rm def}{=}\min\{s,t\}=\int_0^t1_{[0,s]}(u)\,du=\la 1_{0,s]}, 1_{[0,t]}\ra_{L^2([0,1])}
\end{equation}
and let
$$\mch_{00}=\hbox{linear span of $\{K_s\,:\,s\in [0,1]\}$.}$$
It is useful to observe that the derivative ${\dot K}_s$ equals $1$ on $[0,s)$ and $0$ on $(s,1]$:
\begin{equation}\label{E:doths}
{\dot K}_s =1_{[0,s]}\qquad\hbox{everywhere  on $[0,1]$ except at $s=t$.}
\end{equation}
On the vector space $\mch_{00}$ we have the inner product
\begin{equation}\label{E:ipH0}
\la h, k\ra=\int_0^1{\dot h}(t){\dot k}(t)\,dt.
\end{equation}
Then
 \begin{equation}\label{E:hminst}
 \la K_s,K_t\ra  = \min\{s,t\}
 \end{equation}
 for all $s,t\in [0,1]$.  Thus
 $$\la K_s, K_t\ra =K_t(s),$$
 and so, by taking linear combinations,
 \begin{equation}\label{E:repkerhst}
 \la K_s, x\ra = x(s)\qquad\hbox{for all $x\in\mch_{00}$ and $s\in [0,1]$.}
 \end{equation}
 Thus evaluation at $s$ is a linear functional on $\mch_{00}$ that is  given by inner-product against the function $K_s$. The relation  (\ref{E:repkerhst}) is a {\em reproducing kernel} relation, with the function
 $$(s,t)\mapsto \min\{s,t\}$$
 being the reproducing kernel function. Let us also note that point-evaluation is continuous on $\mch_{00}$:
 \begin{equation}\label{E:xtcont}
 |x(s)|=|\la K_s, x\ra| \leq   |\!| K_s|\!||\!|x|\!| =s^{1/2} |\!|x|\!| \leq |\!|x|\!|
 \end{equation}
 for all $x\in\mch_{00}$; in fact,
 \begin{equation}\label{E:supxx}
 |x|_{\sup}\leq  |\!|x|\!|. 
 \end{equation}

 \subsection{Cameron-Martin Hilbert space}\label{ss:cm} We take as the Cameron-Martin Hilbert space the completion of $\mch_{00}$: 
 \begin{equation}\label{E:mch0}\mch_0=\overline{\mch}_{00}.
\end{equation}
We shall now examine the nature of the elements of $\mch_0$.
Suppose $(x_n)$ is a Cauchy sequence in $\mch_0$.  Thus,
$$|\!|{\dot x}_n-{\dot x}_m|\!|_{L^2[0,1]}=|\!|x_n-x_m|\!|\to 0$$
as $n, m\to\infty$, and so there is a function $y\in L^2[0,1]$ for which
$$x_n\to y\quad\hbox{in $L^2[0,1]$.}$$
From the bound (\ref{E:supxx}) it follows that $(x_n)$ is uniformly Cauchy and so converges uniformly to a continuous function $x:[0,1]\to\mbr$; furthermore,
taking $n\to\infty$ in
$$x_n(s)=\la K_s, { \dot x}_n\ra=\int_0^sx_n(s)\,ds,$$
we have
\begin{equation}\label{E:xsys}
x(s)=\la K_s, y\ra=\int_0^s y(u)\,du \qquad\hbox{for all $s\in [0,1]$.}
\end{equation}
Then $x$ is absolutely continuous and, by Lebesgue's theorem, $x$ is  differentiable almost everywhere with $x'=y$ almost everywhere.  The Hilbert space $\mch_{0}$ consists of all absolutely continuous functions $x:[0, 1]\to \mbr$, with initial value $0$ and with $\int_0^1 {\dot x}(t)^2\,dt<\infty$, where the derivative ${\dot x}$ is defined almost everywhere. 

 Finally let us note that the { reproducing kernel} relation (\ref{E:repkerhst}) extends by continuity to 
 \begin{equation}\label{E:repkerhst2}
 \la K_s, x\ra = x(s)\qquad\hbox{for all $x\in\mch_{0}$ and $s\in [0,1]$.}
 \end{equation}

\subsection{The Paley-Wiener map $I$}
 The mapping $I$ is given by the Wiener integral
 \begin{equation}\label{E:Iclass}
 I:\mch_0 \to L^2(\mcc_0,\mu): h\mapsto I(h)=\int_0^1{\dot h}(t)\,dB_t.
 \end{equation}
 This is best understood in our context by examining the values of $I$ on the special elements $K_s\in\mch_{00}$:
 \begin{equation}\label{E:IhsBs}
 I(K_s)=B_s\in\mcc_0^*,
 \end{equation}
 from which we have
 \begin{equation}\label{E:ipIhsht}
 \begin{split}
 \la I(K_s), I(K_t)\ra_{L^2(\mu)} &=\mbe_{\mu}[B_sB_t]\\
 &=\min\{s,t\} \quad\hbox{(by (\ref{E:BsBt}))}\\
 & =\la K_s,K_t\ra_{H}.
 \end{split}
 \end{equation}

  \subsection{The imbedding $i$}
 We also need to understand the imbedding
 $$i:\mch_0\to\mcc_0: h\mapsto i(h)=\int_{\mcc_0}\la h, \omega\ra \omega\,d\mu(\omega).$$
 For any $t\in [0,1]$ we have
 \begin{equation}
 \begin{split}
 h(t) &=\int_0^1  1_{[0,t]}(s){\dot h}(s)\,ds \\
 &=\la  { K}_t, h\ra_E\\
 &=\la  I({  K}_t), i(h)\ra \\
 &\quad\hbox{(evaluating $I({  K}_t)\in \mcc_0^*$ on $i(h)\in\mcc_0$, using (\ref{E:iHBhh_1}))}\\
 & =\la B_t, i(h)\ra \qquad\hbox{ (because $I(K_t)=B_t$ as seen in (\ref{E:IhsBs}))}\\
 &= i(h)(t) \qquad\hbox{ (because $B_t:\mcc_0\to\mbr$ is evaluation at $t$).}
 \end{split}
 \end{equation}
 Thus
 $$i(h)=h,$$
 viewed as an element of $\mcc_0$.  
 
 \subsection{The subspace $\mch_{\mcc_0^*}$ of continuous functionals} Recall that the special function $K_s$ has the property that 
 $$I(K_s)(\omega)=B_s(\omega)=\omega(s)=\int_0^1 \omega\,d\delta_s\qquad\hbox{for all $\omega\in\mcc_0$},$$
 where $\delta_s$ is the delta measure at $s\in [0,1]$.    Taking linear combinations we have, for any $h=\sum_jc_jK_{s_j}\in\mch_{00}$,
\begin{equation}\label{E:Ihomclass}
I(h)(\omega)=\int_0^1\omega \,d\lambda_h\qquad\hbox{for all $\omega\in\mcc_0$},
\end{equation}
 for a  signed   Borel measure $\lambda_h=\sum_jc_j\delta_{s_j}$ on $[0,1]$. A general element of the dual space $\mcc_0^*$ is given by integration against a   Borel signed measure  on $[0,1]$ and corresponds to a function $h\in \mch_0$ for which the derivative ${\dot h}$ is of bounded variation.

  \subsection{Codimension one affine subspaces} 
 
 The simplest type of closed affine subspace of $\mch_0$ is of the form
 $$\{h\in\mch_0\,:\, \la f, h\ra=c\},$$
 for some $f\in \mch_0$ and $c\in\mbr$.  As a special case of interest consider the functon $f=K_1:t\mapsto t$. The reproducing kernel property (\ref{E:repkerhst2}) then implies
 $$\la f,h\ra=\la K_1,h\ra=h(1),$$
 and so the affine subspace above is
 $$  \{h\in\mch_0\,:\,h(1)=c\}.$$
 By (\ref{E:kerfB}) the
  closure in the Banach space $\mcc_0$  of $\ker \la f,\cdot\ra_H$ is $\ker I(f)$:
  \begin{equation}\label{E:closure}
  \{\omega\in\mcc_0\,:\,\omega(1)=0\}=\overline{\{\omega\in\mch_0\,:\,\omega(1)=0\}}
  \end{equation}
  On the left here we have the set of continuous loops based at $0$.

   \subsection{The mapping $I_L$} Consider  a closed affine subspace $L$ of $\mcc_0$. Let $h_L$ be the point on $L$ closest to the origin and 
   $$P_{L_0}:\mch_0\to\mch_0$$
   is the orthogonal projection onto the closed subspace $L_0=L-h_L$. Then, as we have discussed before in the context of (\ref{E:muLdef1}), there is a Gaussian measure $\mu_L$ on $\mcc_0$ for which
   \begin{equation}\label{E:cfmuL}
   \int_{\mcc_0}e^{i  I(h) }\,d\mu_L =e^{i\la h_L, h\ra-\frac{1}{2}|\!|P_{L_0}h|\!|^2} \qquad\hbox{for all $h\in\mch_{\mcc_0^*}$,}   \end{equation}  
   where     \begin{equation}\label{E:Ihomclass2}
I:\mch_0\to L^2(\mcc_0,\mu):h\mapsto \int_0^1{\dot h}(t)\,dB_t
\end{equation}
is the Paley-Wiener map.  Next, for the measure $\mu_L$ we observe  that the mapping
$$\mch_{\mcc_0^*}\to L^2(B,\mu_L): h\mapsto I(h)$$
is continuous linear and extends uniquely to a continuous linear mapping
\begin{equation}\label{E:ILclass}
I_L:\mch_0\to L^2(B,\mu_L).
\end{equation}
We think of $I_L(h)$ as the random variable $I(h)$ {\em conditioned to lie on the affine subspace} $L$. Its distribution is determined by the characteristic function
 \begin{equation}\label{E:cfmuL2}
   \int_{\mcc_0}e^{i  I_L(h) }\,d\mu_L =e^{i\la h_L, h\ra-\frac{1}{2}|\!|P_{L_0}h|\!|^2} \qquad\hbox{for all $h\in\mch_{0}$.}   \end{equation}

 \section{Explicit formulas for transforms}\label{s:expl}

In this section we  work out the Gaussian Radon transforms of  some Brownian functionals. We  use the framework and notation set up in the preceding sections. In particular, we will work with the classical Wiener measure $\mu$ on $\mcc_0$, the Banach space of continuous paths $[0,1]\to\mbr$ initiating at $0$,  the Hilbert space $\mch_0=\overline{\mch}_{00}$ discussed in subsection \ref{ss:cm}, and the continuous linear map 
$$I_L:\mch_0\to L^2(\mcc_0,\mu_L)$$
for closed affine subspaces $L\subset \mcc_0$. 

We will determine $G\phi$ for $\phi$   a  linear combination of  multiple It\^o integrals 
$$U(g^{\otimes n})= \int_0^1\ldots\int_0^1 g(s_1)\ldots g(s_n)\,dB_{s_1}\ldots dB_{s_n},$$
where $g\in L^2[0,1]$ and $U$ is the Wiener chaos isomorphism
\begin{equation}
U: \mcf_s(H)\to L^2(\mch_0) 
\end{equation}
discussed in (\ref{E:fockspaceu}).

\subsection{An affine subspace of pinched paths} Let us consider a concrete and simple example. Let us take for $L$ the closed affine subspace
\begin{equation}\label{E:LTc1}
L_T(c)= \{h\in\mch_0\,: h(T)=c\}=\{h\in\mch_0\,: \la K_T, h\ra=c\},
 \end{equation}
where $T$ is a fixed time in $(0,1]$ and $c\in\mbr$, and, as usual,
$$K_T:[0,1]\to\mbr: s\mapsto \min\{T,s\}.$$
Geometrically, a point $h$ in the affine subspace $L_T(c)$ is a path  initiating at $0$  and passing through the point $c$ at time $T$.  If $T=1$ this would be a {\em bridge} running from $0$ to $c$. 

We can display $L_T(c)$ also as
 \begin{equation}\label{E:LTc2}
 L_T(c)= \frac{c}{ {T}}K_T + K_T^{\perp},
 \end{equation}
 where $T$ in the first term on the right side arises as $|\!|K_T|\!|^2$ and
the point on $L_T(c)$ closest to the origin is
 \begin{equation}\label{E:LTc3}
 h_{L_T(c)}= \frac{c}{ {T}}K_T.
 \end{equation}
 The orthogonal projection $P_{L_T(c)}:\mch_0\to\mch_0$ onto the closed subspace $K_T^{\perp}$ is given by
\begin{equation}\label{E:PLTch}
P_{L_T(c)}(h) = h -\frac{1}{T}\la h,K_T\ra K_T.
\end{equation}
\subsection{The bridge process}
Let us recall the linear mapping 
$$I_{L_T(c)}:\mch_0\to L^2(\mcc_0,\mu_L)$$
given in
(\ref{E:defILh}).
The distribution of the random variable 
$$I_{L_T(c)}(h)\in L^2(\mcc_0,\mu_L)$$
is determined through its characteristic function, which is given by (\ref{E:cfmuL2}) as
 \begin{equation}\label{E:cfmuL3}
   \int_{\mcc_0}e^{i uI_{L_T(c)}(h) }\,d\mu_{L_T(c)} =e^{iu\frac{c}{T}\la K_T, h\ra- \frac{u^2}{2}|\!| h -\frac{1}{T}\la h,K_T\ra K_T|\!|^2}    \end{equation}  
 for all $h\in\mch_{0}$ and $u\in\mbr$. The variance term in the exponent works out to
\begin{equation}\label{E:varh}
{\rm var}\Bigl(I_{L_T(c)}(h) \Bigr)= |\!|h|\!|^2-\frac{1}{T}\la h,K_T\ra^2.
\end{equation}
 Recalling from (\ref{E:Iclass}) that 
\begin{equation}\label{E:IhinthBdot}
I(h)= \int_0^1{\dot h}(t)\,dB_t,
\end{equation}
 let us use the notation
 \begin{equation}\label{E:ILh}
 \ I_{L_T(c)}(h)  \stackrel{\rm def}{=} \int_0^1{\dot h}(t)\,dB_t \qquad\hbox{for $h\in\mch_0$.}
 \end{equation}
 There is a slight conflict between this usage and  (\ref{E:Iclass}) because $I_{L_T(c)}(h)$ is an element of $L^2(\mcc_0,\mu_{L_T(c)})$ whereas $I(h)$ is an element of the Wiener $L^2$-space $L^2(\mcc_0,\mu)$. But, as noted in (\ref{E:ImuL}),  the mapping  $I_{L_T(c)}$ is defined to agree with $I$
 on the elements $h\in\mch_{00}$. Furthermore, $I(h)$ is a {\em continuous} function on $\mcc_0$ when $h\in\mch_{00}$ (which is contained in $H_{B^*}$ in this context) and so is uniquely specified not simply almost-everywhere but  at all points of $\mcc_0$. 
 
 Let $f={\dot h}$, where $h\in \mch_0$; the characteristic function given in (\ref{E:cfmuL3})  and the variance formula (\ref{E:varh}) show  that $ \int_0^1f(t)\,dB_t$ is Gaussian with mean and variance
 \begin{equation}\label{E:meanvarILh}
 \begin{split}
 \mbe_{\mu_{L_T(c)}}\Bigl( \int_0^1f(t)\,dB_t\Bigr)&= \frac{c}{T}\int_0^Tf(t)\,dt\\
{\rm var}_{\mu_{L_T(c)}}\Bigl( \int_0^1f(t)\,dB_t\Bigr)&= \int_0^1f(t)^2\,dt -\frac{1}{T}\Bigl(\int_0^T f(t)\,dt\Bigr)^2.\end{split}
 \end{equation}
 If we take for $h$ the function $K_t$,  then $f=1_{[0,t]}$ and so
    $$I(K_t)=\int_0^1{\dot K}_t(s)\,dB_s=\int_0^tdB_s=B_t.$$ 
Then $I_{L_T(c)}(K_t)$ is Gaussian with mean and variance given by
 \begin{equation}\label{E:EvarLTcHt}
 \begin{split}
 \mbe_{\mu_{L_T(c)}}\bigl(B_t\bigr) &= \frac{c}{T}\la K_T, K_t\ra=\frac{c}{T}\min\{T,t\}\\
 {\rm var}_{\mu_{L_T(c)}}\bigl(B_t\bigr) &= t-\frac{1}{T}\min\{t,T\}^2=T\Bigl(\frac{t}{T}-\min\{1,t/T\}^2\Bigr).
 \end{split}
 \end{equation}
 The expectation term is the Gaussian Radon transform; thus the value of the Gaussian Radon transform $GB_t$ on the affine subspace $L_T(c)$
 is
 \begin{equation}\label{E:LTch}
 GB_t|_{[B_T=c]}= \frac{c}{T}\min\{T,t\},
 \end{equation}
 where we have sacrificed some notational accuracy for notational intuition in writing $[B_T=c]$ to denote the affine subspace $L_T(c)$.
 
Using the covariance formula (\ref{E:covILH12}) we have
 \begin{equation}\label{E:covILH12class}
 \begin{split}
{\rm Cov}_{\mu_{L_T(c)}}\bigl(B_t, B_s\bigr) &= \la P_{L_0}K_t, P_{L_0}K_s\ra\\
&=\left\la K_t -\frac{1}{T}\la K_t,K_T\ra K_T, K_s -\frac{1}{T}\la K_s,K_T\ra K_T \right\ra\\
&=\la K_t,K_s\ra-\frac{1}{T}\la K_t,K_T\ra\la K_s, K_T\ra \\
&=\min\{s,t\} -\frac{1}{T}\min\{t,T\}\min\{s,T\}\\
&=T\Bigl(\min\left\{\frac{s}{T}, \frac{t}{T}\right\} - \min\left\{1, \frac{t}{T}\right\} \min\left\{1, \frac{s}{T}\right\}\Bigr).
\end{split}
\end{equation}
This covariance structure identifies the process
$$t\mapsto I_{L_T(c)}(K_t)$$
as a {\em Brownian bridge process} for $t\in [0,T]$ and a Brownian motion for $t>T$ starting at $c$ when $t=T$. 

\subsection{Multiple bridges} Instead of conditioning the path to pass through a specified point at a particular time $T$ we can   condition the process to pass through a sequence of specific points at specified times $T_1<\ldots<T_m$.  For this the relevant affine subspace of $\mch_0$ is
\begin{equation}\label{E:LTmulti}
L_T(c)=\{h\in\mch_0\,:\, h(T_1)=c_1,\ldots, h(T_m)=c_m\ra
\end{equation}
where now $T=(T_1,\ldots, T_m)\in [0,1]^m$ with $0<T_1<\ldots <T_m$ and $c$ is the point $(c_1,\ldots, c_m)\in\mbr^m$. The special point $h_{L_T(c)}$ on $L_T(c)$ closest to the origin is given by the formula in Proposition \ref{P:hLforcodimfin}.   The affine subspace $L_T(c)$ is the level set  
$$L_T(c) =F^{-1}(c),$$
where
\begin{equation}\label{E:Fbridge}
F:\mch_0\to\mbr^m: h\mapsto \Bigl(\la K_{T_1}, h\ra,\ldots, \la K_{T_m},h\ra\Bigr).
\end{equation}
By Proposition \ref{P:hLforcodimfin} the matrix $FF^*$ has $(j,k)$-the entry given by
\begin{equation}\label{E:FFstarcw}
(FF^*)_{jk}=\la K_{T_k}, K_{T_j}\ra=\min\{T_j, T_k\}.
\end{equation}
Using this and the expression for $GI(h)$ obtained in (\ref{E:GIhhLc}) we have
\begin{equation}\label{E:GIHcw}
GI(h)|_{L_T(c)}=\sum_{j=1}^m[(FF^*)^{-1}]_{jk}c_k\la K_{T_j},h\ra.
\end{equation}
 A more transparent formula is obtained if we use a different representation for the affine subspace $L_T(c)$. To this end let
 $$ w_1=\frac{K_{T_1}}{\sqrt{T_1}},\qquad  w_2=\frac{K_{T_2}-K_{T_1}}{\sqrt{T_2-T_1}} ,\quad\ldots, \quad w_m= \frac{K_{T_m}-K_{T_{m-1}}}{\sqrt{T_m-T_{m-1}}} . $$
 Then
 $$L_T(c)=\{h\in\mch_0: \la w_1,h\ra=b_1,\ldots, \la w_m,h\ra=b_m\}$$
 where
 $$b_1=\frac{c_1}{\sqrt{T_1}},\qquad b_2=\frac{c_2-c_1}{\sqrt{T_2-T_1}} ,\ldots, b_m= \frac{c_{m}-c_{{m-1}}}{\sqrt{T_m-T_{m-1}}}.$$
 The advantage of this representation for $L_T(c)$ is that the vectors $w_1,\ldots, w_m$ are orthonormal.  Then using the formula for $GI(h)$ given in (\ref{E:GIhhLcorth}) we have
\begin{equation}\label{E:GIhLTc2}
\begin{split}
GI(h)|_{L_T(c)} &=\sum_{k=1}^m b_k\la w_k,h\ra\\
&=\frac{c_1}{T_1}\int_0^{T_1}{\dot h}(t)\,dt+\sum_{k=2}^m \frac{c_k-c_{k-1}}{T_k-T_{k-1}}\int_{T_{k-1}}^{T_k}{\dot h}(t)\,dt.
\end{split}
\end{equation}
Expressing this in different notation, and writing $f$ for $\dot h$, we have
\begin{equation}\label{E:GIhLTc3}
\mbe_{\mu_{L_T(c)}}\Bigl(\int_0^1f(t)\,dB_t\Bigr)   =\frac{c_1}{T_1}\int_0^{T_1}f(t)\,dt+\sum_{k=2}^m \frac{c_k-c_{k-1}}{T_k-T_{k-1}}\int_{T_{k-1}}^{T_k}{f}(t)\,dt.
\end{equation}
 This generalizes the first equation in (\ref{E:EvarLTcHt}). As a check we observe that when $f=1_{[0,T_1]}$ the formula (\ref{E:GIhLTc3}) yields
 $$ \mbe_{\mu_{L_T(c)}}\Bigl(B_{T_1}\Bigr)=c_1,$$
 and, more generally, taking $f=1_{[0,T_j]}$, we obtain
  $$ \mbe_{\mu_{L_T(c)}}\Bigl(B_{T_k}\Bigr)=c_1+\sum_{k=2}^j \frac{c_k-c_{k-1}}{T_k-T_{k-1}}(T_k-T_{k-1})=c_j,$$
  which is consistent with the process $t\mapsto B_t$ passing through the point  $c_k$ at time $T_k$, for each $k\in\{1,\ldots, m\}$, almost surely with respect to the measure $\mu_{L_T(c)}$.
 
 \subsection{Bridges in higher dimension} For this let us simply note that for an $\mbr^n$-valued Brownian motion process the Cameron-Martin Hilbert space is $\mch_0^n$, where $\mch_0$ is the Hilbert space for the one-dimensional process. A  bridge process, or more precisely a process that is required to pass through a point $(p_1,\ldots, p_n)\in\mbr^n$ at a specified time $T$, corresponds then to a closed affine subspace of $\mch_0^n$ specified by level sets of the  $n$ orthogonal vectors in $\mch_0$ given by
 $$v_1=(K_T,0,\ldots, 0), \ldots, v_n=(0,\ldots, 0, K_{T})\in\mch_0^n.$$
 We can obtain Gaussian Radon transforms of functions of the form $I(h)$ by applying Proposition \ref{P:hLforcodimfin} in a manner similar to what we did in the preceding subsection. We shall not explore this further.

\subsection{Multiple stochastic integrals} We refer to the monograph of Kuo \cite[sec. 9.6]{Ku2} for a development of It\^o's theory of multiple stochastic integrals.
For any symmetric function  $F\in L^2[0,1]^n$,   the multiple It\^o stochastic integral is given by $n!$ times the iterated integral; thus:
 \begin{equation}\label{E:multito}
 \begin{split}
J_n(F)&\stackrel{\rm def}{=} \int_{[0,1]^n} F(t_1,\ldots, t_n)\,dB_{t_n} \ldots dB_{t_1}\\
&
\stackrel{\rm def}{=} n!\int_0^1\ldots\int_0^{t_{n-2}}\left[\int_0^{t_{n-1}}F(t_1,\ldots, t_n)\,dB_{t_n}\right]\ldots dB_{t_1}.
\end{split}
\end{equation}
For any $f\in L^2[0,1]$,  the integral of $f^{\otimes n}$ can be expressed as a Hermite polynomial:
\begin{equation}\label{E:Jnfint}
J_n(f^{\otimes n})=H_n\left(\int_0^1f(s)\,dB_s;\, |\!|f|\!|_{L^2[0,1]}^2\right) 
\end{equation}
where the Hermite polynomial
$$H_n(x; u^2)$$
is specified through the generating formula
\begin{equation}\label{E:hermitegen}
e^{tx-\frac{u^2}{2}t^2}=\sum_{n=0}^\infty \frac{t^n}{n!}H_n(x; u^2).
\end{equation}
(For a proof of (\ref{E:multito}) and more details  we refer  to \cite{Ku2}.)
The left side remains unchanged under the following transformation:
\begin{equation}
\Bigl(t,x,u^2\Bigr)\mapsto \Bigl({\lambda}^{-1}t, {\lambda} x, {\lambda} u^2\Bigr) 
\end{equation}
for any $\lambda\in\mbr$. 
As a result we have
\begin{equation}\label{E:scale}
H_n\Bigl({\lambda}x; {\lambda}^2u^2\Bigr)={\lambda}^n H_n(x; u^2).
\end{equation}

Using the generating formula  (\ref{E:hermitegen}) we have 
\begin{equation}\label{E:eIHHn}
e^{I(h)-\frac{1}{2}|\!|h|\!|_0^2}  = \sum_{n=0}^{\infty} \frac{1}{n!} H_n\Bigl(I(h); |\!|h|\!|_0^2\bigr).
\end{equation}
Writing $f={\dot h}$, which is in $L^2[0,1]$ and using
$$I(h)=\int_0^1{\dot h}(t)\,dB_t=\int_0^1f(t)\,dB_t$$
and
$$|\!|h|\!|_0^2=|\!|{\dot h}|\!|_{L^2[0,1]}^2=|\!|f|\!|_{L^2[0,1]}^2,$$
we have from (\ref{E:eIHHn}) and (\ref{E:Jnfint}):
\begin{equation}\label{E:eIHHn2}
e^{I(h)-\frac{1}{2}|\!|h|\!|_0^2} =\sum_{n=0}^\infty\frac{1}{n!}J_n({\dot h}^{\otimes n}).
\end{equation}
We have noted that for $h\in\mch_{00}$  the Wiener integral $\int_0^1{\dot h}(s)\,dB_s$ is a continuous function on $\mcc_0$, and hence so is $J_n({\dot h}^{\otimes n})$. This connects up with 
the Fock space isomorphism
\begin{equation}
U: \mcfs(\mch_0)_c\to L^2(\mcc_0;\mu): \Exp(h)\mapsto e^{I(h)-\frac{1}{2}|\!|h|\!|_0^2}\end{equation}
for $h\in\mch_0$.  

\subsection{Gaussian Radon transform of multiple It\^o integrals}
The Gaussian Radon transform of $J_n({\dot h}^{\otimes n})$ at the closed affine subspace $L_T(c)$ is
\begin{equation}\label{E:EJnLTc}
\begin{split}
GJ_n({\dot h}^{\otimes n})\Big|_{L_T(c)}&=\mbe_{\mu_{L_T(c)}}\Bigl[
H_n\left(\int_0^1f(s)\,dB_s; |\!|f|\!|_{L^2[0,1]}^2\right)
\Bigr]\\ 
&=\mbe\left[H_n\Bigl(I_{L_T(c)}({ h}); |\!|{\dot h}|\!|_{L^2[0,1]}^2 \Bigr) \right].
\end{split}
\end{equation}

Let $X$ be any Gaussian variable.  Then from the generating formula (\ref{E:hermitegen}) for the Hermite polynomials we have
\begin{equation}\label{E:etXHn}
\mbe\Bigl(e^{tX-\frac{u^2}{2}t^2}\Bigr) =\sum_{n=0}^\infty\frac{t^n}{n!}\mbe[H_n(X; u^2)]
\end{equation}
for all $t, u\in\mbr$.
On the other hand
from     the Gaussian formula
\begin{equation}\label{E:Gaussfunda2}
\mbe\bigl[e^{tX}\bigr]=e^{\mbe[X]t+\frac{1}{2}{\rm var}(X)t^2}
\end{equation}
we have 
\begin{equation}\label{E:Gaussfunda}
\mbe\bigl[e^{tX-\frac{u^2}{2}t^2}\bigr]=e^{\mbe[X]t-\frac{1}{2}\bigl(u^2-{\rm var}(X) \bigr)t^2}.
\end{equation}
Writing
$$v^2=u^2-{\rm var}(X)$$
we then have
\begin{equation}\label{E:Gaussfunda3}
\mbe\bigl[e^{tX-\frac{u^2}{2}t^2}\bigr]=\sum_{n=0}^\infty\frac{t^n}{n!}H_n\Bigl(\mbe(X); v^2\Bigr).
\end{equation}
Here we assume that $u^2\geq {\rm var}(X)$.
Comparing with (\ref{E:etXHn}), we have the following  shift-of-variance relation
\begin{equation}\label{E:shiftvar}
\mbe\bigl(H_n(X; u^2)\bigr)= H_n\Bigl(\mbe(X); v^2\Bigr).
\end{equation}
Now let us apply this to the Gaussian Radon transform of the iterated integral of $J_n({\dot h}^{\otimes n})$, for $h\in\mch_{00}$, given in (\ref{E:EJnLTc}):
\begin{equation}\label{E:EJnLTc2}
\begin{split}
GJ_n({\dot h}^{\otimes n})\Big|_{L_T(c)}&= \mbe\left[H_n\Bigl(I_{L_T(c)}({  h}); |\!|{\dot h}|\!|_{L^2[0,1]}^2 \Bigr) \right]\\
&= H_n\Bigl(\mbe\bigl(I_{L_T(c)}({  h})\bigr);  |\!|{\dot h}|\!|_{L^2[0,1]}^2-{\rm var}\bigl(I_{L_T(c)}({  h})\bigr)\Bigr) \\
&=  H_n\left(\frac{c}{T}\int_0^T{\dot h}(t)\,dt; \frac{1}{T}\left(\int_0^T{\dot h}(t)\,dt\right)^2\right)\\
&\qquad\hbox{(using (\ref{E:meanvarILh})) }\\
&=H_n\left(\frac{c}{T}h(T);\frac{1}{T}h(T)^2\right)\\
&=h(T)^nH_n\left(\frac{c}{T};\frac{1}{T}\right)
\end{split}
\end{equation}
where in the last line we used the scaling relation (\ref{E:scale}). 

Writing
$$f={\dot h},$$
we can rewrite this as
\begin{equation}\label{E:Giterated}
GJ_n(f^{\otimes n})|_{L_T(c)}=\left(\int_0^Tf(t)\,dt\right)^nH_n\left(\frac{c}{T};\frac{1}{T}\right).
\end{equation}
This is the Gaussian Radon transform of the mutliple It\^o integral
$$n!\int_0^1\ldots\int_0^{t_{n-2}}\left[\int_0^{t_{n-1}}f(t_1)\ldots f(t_n)\,dB_{t_n}\right]\ldots dB_{t_1}$$
evaluated at the affine subspace given by
$$L_T(c)=\{h\in \mch_0: h(T)=c\},$$
corresponding geometrically to paths that start  at the point $0$ and pass through $c$ at time $T$. In the special case of $T=1$ and $f=1$, the multiple integral $J_n(f^{\otimes n})$ is, by (\ref{E:Jnfint}), equal to
$$J_n(1^{\otimes n})=H_n(B_1; 1).$$
Now $B_1=c$ almost surely with respect to the Brownian bridge measure $\mu_{L_1(c)}$ and so $J_n(1^{\otimes n})=H_n(c)$, the traditional Hermite polynomial evaluate at $c$.  This coincides with the value given by the right hand side of (\ref{E:Giterated}).

We observe that the integral on the right hand side of (\ref{E:Giterated}) is
$$\int_{[0,T]^n}f(t_1)\ldots f(t_n)\,dt_1\ldots dt_n.$$
By the algebraic process of polarization we can express a symmetric tensor product
\begin{equation}
f_1{\hat\otimes}\ldots{\hat\otimes}f_n= \frac{1}{n!}\sum_{\sigma\in S_n}f_{\sigma^{-1}(1)}\otimes\ldots \otimes f_{\sigma^{-1}(n)}
\end{equation}
as a linear combination of terms of the form $g^{\otimes n}$, where the elements $g$ are linear combinations of $f_1,\ldots, f_n$. Then from (\ref{E:Giterated})  we obtain the more general form
\begin{equation}\label{E:Giterated2}
GJ_n(f_1{\hat\otimes}\ldots{\hat\otimes}f_n)|_{L_T(c)}=\left(\int_{[0,T]^n}f_1(t_1)\ldots f_n(t_n)\,dt_1\ldots dt_n\right) H_n\left(\frac{c}{T};\frac{1}{T}\right).
\end{equation}
Here each $f_j$ is of the form ${\dot h}_j$, where $h_j\in \mch_{\mcc_0^*}$, and in particular could be in the subspace $\mch_{00}\subset\mch_0$.  Writing $F$ for $ f_1{\hat\otimes}\ldots{\hat\otimes}f_n$, viewed as an element of $L^2([0,1]^n)$, or taking $F$ to be a linear combination of such functions, we then have
\begin{equation}\label{E:Giterated3}
GJ_n(F)|_{L_T(c)}=\left(\int_{[0,T]^n}F(t_1,\ldots, t_n)\,dt_1\ldots dt_n\right) H_n\left(\frac{c}{T};\frac{1}{T}\right).
\end{equation}
Here on the left  we have the Gaussian Radon transform of the multiple It\^o stochastic integral $\int_{[0,1]^n}F(t_1,\ldots, t_n)\,dB_{t_n}\ldots dB_{t_1}$ as defined in  (\ref{E:multito}), evaluated on the affine subspace $L_T(c)$ corresponding to Brownian paths satisfying $B(T)=c$.

\section{Concluding Remarks}

 In this paper we have studied Gaussian measure conditioned to be supported on closed affine subspaces of a Hilbert space. We proved a Fock space decomposition for such measures. We examined the Gaussian Radon transform, with special attention to finite-codimension affine subspaces. We investigated the transform and the conditioned measure for the case of the classical Wiener space, observing how Brownian bridges can be understood in terms of affine subspaces of finite codimension. We computed Gaussian Radon transforms of multiple It\^o stochastic integrals explicitly in terms of Hermite polynomials.

{\bf Acknowledgments}. This work is part of a research project covered by  NSA grant
  H98230-13-1-0210. 
  	
 \bibliographystyle{amsplain}

\end{document}